\documentclass[twoside]{IEEEtran}

\usepackage{atbegshi,picture,ifpdf}
\usepackage[noadjust]{cite}
\usepackage[]{xcolor}
\usepackage[pdftex]{graphicx}
\usepackage[cmex10]{amsmath}
\usepackage[caption=false,font=footnotesize]{subfig}
\usepackage{stfloats}

\usepackage[english]{babel}
\usepackage[utf8x]{inputenc}
\usepackage[T1]{fontenc}
\usepackage[scaled=0.91]{helvet}
\usepackage[cmintegrals]{newtxmath}
\usepackage{acronym}
\usepackage{bm}

\IEEEeqnarraydefcolsep{0}{\leftmargini}

\newcommand*\mat[1]{\textsf{\textit{\textbf#1}}}

\renewcommand*\vec[1]{\small\textsf{\textit{\textbf#1}}}
\newcommand*\vecelem[1]{\small\textsf{\textit{#1}}}

\makeatletter
\def\ps@IEEEtitlepagestyle{%
  \def\@oddfoot{\mycopyrightnotice}%
  \def\@oddhead{\hbox{}\@IEEEheaderstyle\leftmark\hfil\thepage}\relax
  \def\@evenhead{\@IEEEheaderstyle\thepage\hfil\leftmark\hbox{}}\relax
  \def\@evenfoot{}%
}
\def\mycopyrightnotice{%
  \begin{minipage}{\textwidth}
  \centering \scriptsize
  Copyright~\copyright~2018 IEEE. Personal use of this material is permitted. Permission from IEEE must be obtained for all other uses, in any current or future media, including\\reprinting/republishing this material for advertising or promotional purposes, creating new collective works, for resale or redistribution to servers or lists, or reuse of any copyrighted component of this work in other works by sending a request to pubs-permissions@ieee.org.
  \end{minipage}
}
\makeatother

\AtBeginShipout{\AtBeginShipoutUpperLeft{%
  \put(\dimexpr\paperwidth-1cm\relax,-.55cm){\makebox[0pt][r]{
  \begin{minipage}{\textwidth}
  \centering \scriptsize 
  This is the author's version of an article that has been published in this journal. 
  Changes were made to this version by the publisher prior to publication. 
  \\[0.3ex]
  The final version of record is available at\qquad http://dx.doi.org/10.1109/TAP.2018.2866578
  \end{minipage}
  }}}%
}

\begin{document}
\title{An Accurate Low-Order Discretization Scheme for the Identity Operator in the Magnetic Field and Combined Field Integral Equations}

\author{Jonas Kornprobst, \IEEEmembership{Student Member, IEEE}, and Thomas F. Eibert, \IEEEmembership{Senior Member, IEEE}%
\thanks{Manuscript received December 4, 2017; revised July 27, 2018; accepted August 8, 2018;
date of this version August 8, 2018. \emph{(Corresponding author: Jonas Kornprobst.)}}%
\thanks{J. Kornprobst and T. F. Eibert are with the Chair of High-Frequency Engineering, Department of Electrical and Computer Engineering, Technical University  of  Munich,  80290  Munich, Germany (e-mail:  j.kornprobst@tum.de, hft@ei.tum.de).}%
\thanks{Color versions of one or more of the figures in this paper are available online at http://ieeexplore.ieee.org.}%
\thanks{Digital Object Identifier 10.1109/TAP.2018.286657}%
}

\markboth{IEEE Transactions on Antennas and Propagation}%
{Kornprobst and Eibert: An Accurate Low-Order Discretization Scheme for the Identity Operator in the MFIE and CFIE}

\maketitle

\begin{abstract}
A new low-order discretization scheme for the identity operator in the magnetic field integral equation (MFIE) is discussed. 
Its concept is derived from the weak-form representation of combined sources which are discretized with Rao-Wilton-Glisson (RWG) functions. 
The resulting MFIE overcomes the accuracy problem of the classical MFIE while it maintains fast iterative solver convergence. 
The improvement in accuracy is verified with a mesh refinement analysis and with near- and far-field scattering results. 
Furthermore, simulation results for a combined field integral equation (CFIE) involving the new MFIE show that this CFIE is interior-resonance free and well-conditioned like the classical CFIE, but also accurate as the EFIE. 
\end{abstract}


\section{Introduction} %
\IEEEPARstart{E}{lectromagnetic} radiation and scattering problems with perfect electrically conducting (PEC) objects are commonly treated by the method of moments (MoM) in conjunction with boundary integral equations due to their good accuracy and low complexity, i.e.\ surface discretization only. 
The electric field integral equation (EFIE) works with the tangential electric fields on the surface of a scatterer and can be applied for open and closed structures, whereas the magnetic field integral equation (MFIE) for the tangential magnetic field on the surface is suitable for closed bodies only~\cite{peterson1998computational}. 
If applied to closed bodies, the solution to both EFIE and MFIE with the most general form of electric and magnetic surface currents is not unique without side constraints. 
This issue is usually overcome by the choice of specific equivalent sources, the so-called Love-currents~\cite{Schelkunoff1936,Love1901}. 

In addition, there is a second uniqueness problem: 
Both EFIE and MFIE suffer from interior resonances of the corresponding PEC or perfect magnetically conducting (PMC) cavities, respectively~\cite{Maue_1949,Mei_1963,Andreasen_1964,Peterson1990}. 
This can lead to a null space in the MoM equation system and, thus, very poor iterative solver convergence and solution accuracy. 
Simply adding EFIE and MFIE to the so-called combined field integral equation (CFIE) is the most common technique to solve this issue~\cite{mautz1979h}. 
An important benefit of utilizing the CFIE is its improved iterative solver convergence as compared to pure first kind EFIE solutions, as it contains a second kind Fredholm integral equation in the form of the MFIE.

Surface modeling is commonly working with a triangular mesh. 
The simplest divergence-conforming basis functions for the representation of surface current densities on such a mesh are the Rao-Wilton-Glisson (RWG) functions~\cite{Rao_1982}. 
The evaluation of the fields on the surface is most advantageously performed with testing functions from the dual space of the range of the considered operator acting on the surface current densities~\cite{Yla-OijalaKiminkiMarkkanenEtAl2013,Yla-OijalaMarkkanenJarvenpaaEtAl2014}. 
However, especially with the low-order RWG functions, it is not always feasible to follow these rules. 
Also, it was found that the MFIE suffers from inaccuracy problems, in particular due to the involved identity operator~\cite{ErguelGuerel2009}. 

In~\cite{Rius_2001,Ergul_2004}, the {MFIE} was shown to exhibit bad accuracy with RWG basis function discretization, especially for sharp edges and electrically small features. 
It is clear that these errors also contaminate the results of the CFIE~\cite{Gurel_2009}. 
In~\cite{ubeda2005mfie,Ergul_2006}, curl-conforming basis functions such as the $\bm n \times \textrm{RWG}$ functions were studied to improve the accuracy of the MFIE, since the integral operator for electric surface currents in the MFIE does not really require div-conforming basis functions~\cite{Peterson}. 
However, in the more general dielectric case, electric fields must also be computed from these functions with the need to appropriately handle additional hyper-singular line charges along the triangle edges~\cite{Rao1990}. 
The same problem arises when monopolar RWG functions are considered~\cite{ubeda2005monopolar,Ubeda2006monopolar,Zhang2010Monopolar}. 
Nevertheless, these approaches do not offer solutions with satisfying accuracy. 
Full first-order discretization with twice as many unkowns can provide higher accuracy~\cite{Ergul2006linearlinear,Ergul_2007,Ismatullah_2009_MFIE}.

Over the years, many different types of basis and testing functions as well as various discretization schemes have been employed  to analyze and enhance the accuracy of the {MFIE}~\cite{Davis2004higherorder,Davis2005error,Ylae-OijalaTaskinenJaervenpaeae2005,Warnick2007regularization,Ylae-OijalaTaskinenJaervenpaeae2008,Peterson2008,Cools_2010,Cools_2011,Ylae-OijalaKiminkiCoolsEtAl2012,YanJinNie2011,YanJinNie2013,Li2012,Pan2014,Karaos20016,Huang2016solidangle}.
The different combinations of basis and testing functions influence the accuracy of the {MFIE} mainly due to their different representations of the identity operator. 
With low-order basis functions, well-conditioned MFIE formulations might not be able to achieve proper testing of the identity operator~\cite{Yla-OijalaKiminkiMarkkanenEtAl2013,Yla-OijalaMarkkanenJarvenpaaEtAl2014,ErguelGuerel2009}. 
The only exception are the quasi-curl- and div-conforming Buffa-Christiansen (BC) functions, which have been investigated with respect to their suitability for the discretization and testing of the MFIE and improved results for this mixed discretization have been reported~\cite{chen1990,buffa2007dual,Cools_2010,Cools_2011,Ylae-OijalaKiminkiCoolsEtAl2012,Yan2011identity,YanJinNie2011,YanJinNie2013}. 
However, RWG functions are still advantageous with regard to applicability and simulation effort.

As an alternative to the CFIE, the combined source integral equation (CSIE) is also known to be stable at interior resonances~\cite{Bolomey_1973,Mautz1979,Rogers1985,Lee2005,Morita1990IE}.  
In the mathematic community, the combined-source (CS) condition is formulated as a relation between the potential and its normal derivative, known as the Brakhage-Werner trick~\cite{brakhage1965dirichletsche,buffa2005regularized,darbas2006generalized,steinbach2009modified,melenk2012mapping}.
The implementation of the CSIE has been studied with RWG functions for the electric currents and BC functions for magnetic currents~\cite{Yla_Oijala_2012,Ylae-OijalaKiminkiJaervenpaeae2013}. 
Alternatively, RWG functions for both types of equivalent sources are possible~\cite{kornprobst_accurate_2017,kornprobst_accurate_stable_2017,Kornprobst_2017,kornprobst_RWG_2017}. 
The CS condition has been discretized in weak-form in all these works, since weak-form implementations are widely known to be accurate for the similar Leontovich impedance boundary condition~\cite{Ismat2009MoM, YanJin2013}. 
The weak-form conditions are advantageous, because the strong form leads to problems with hyper-singular line charges in the very near field. While this is not really problematic in larger distances~\cite{EibertVojvodicHansen2016}, the weak-form is still more accurate~\cite{Eibert2017}.
 
\begin{figure*}[t]
\centering
 \begin{minipage}[b]{0.245\linewidth}
 \centering%
  \subfloat[]{%
 \includegraphics[]{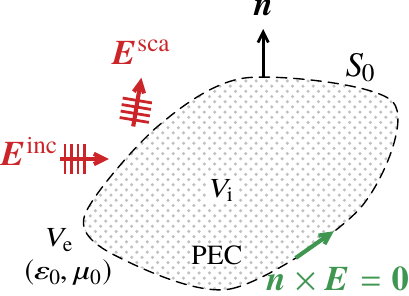}%
   }%
 \end{minipage}%
 \hfill
 \begin{minipage}[b]{0.245\linewidth}
 \centering%
  \subfloat[]{%
 \includegraphics[]{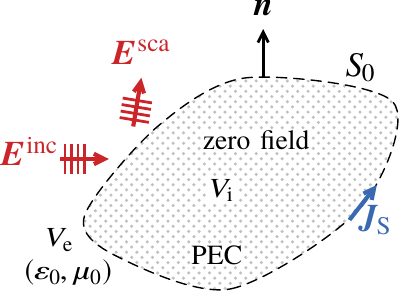}%
   }%
 \end{minipage}%
 \hfill
 \begin{minipage}[b]{0.245\linewidth}
 \centering%
  \subfloat[]{%
 \includegraphics[]{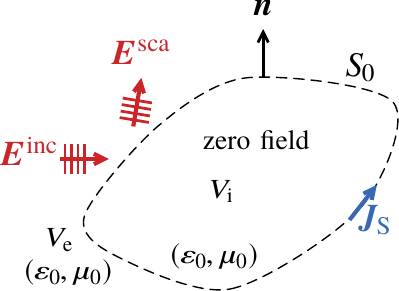}%
   }%
 \end{minipage}%
 \hfill
 \begin{minipage}[b]{0.245\linewidth}
 \centering%
  \subfloat[]{%
 \includegraphics[]{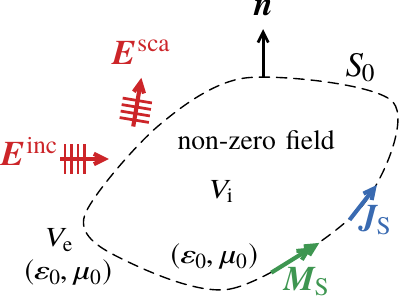}%
   }%
 \end{minipage}
 \caption{The scattering problem with incident wave~$\bm E^\mathrm{inc}$ and scattered field~$\bm E^\mathrm{sca}$ for a PEC scatterer. (a) Scattering object inside the source domain $V_\mathrm i$ embedded in free space. (b) Scatterer with radiating Love currents~$\bm J_\mathrm S$ and zero field in $V_\mathrm i$. (c) Source domain replaced by free space, Love currents radiating in free space. (d) Other choices than the Love currents result in non-zero fields inside~$V_\mathrm i$.\label{fig:scatterer}}
 \vspace*{-0.3cm}
\end{figure*}%
 
Even though there is an increase in accuracy, the CSIE has several disadvantages compared to the CFIE: 
The number of unknowns is doubled and the MoM matrix parts for electric and magnetic currents have to be stored separately. 
This increases the time per matrix vector product, even though the memory consumption is not affected (both matrices are symmetric and can be stored efficiently). 
A further disadvantage is that the Love-currents as utilized in the CFIE are more flexible in handling different kinds of problems, e.g. in hybrid finite-element boundary-integral methods. 
Thus, it is still desirable to work with the CFIE in some cases. 

In this paper, we derive a novel accurate weak-form discretization scheme for the identity operator, which is fully based on RWG functions. 
This scheme is motivated from the previously studied weak-form CS condition and utilizes the corresponding weak-form vector rotations.
The new identity operator discretization leads to an accurate discretization of the MFIE with RWG functions only, and subsequently also to an accurate discretization of the CFIE.

The paper is structured as follows. 
First, we review the integral equations and their discretization for the PEC scenario, with focus on the discretization of the recently proposed weak-form CSIE  with RWG functions. 
Based on the corresponding weak-form rotation, we introduce a basis transformation scheme for the classically discretized MFIE. 
Different transformation schemes for the basis and testing functions are discussed, and a suitable discretization of the identity operator is found to solve the accuracy problem of the MFIE. 
The mesh refinement behavior and solution accuracy of the proposed MFIE formulation are analyzed in near-field and far-field. 
Finally, the new MFIE is applied within the CFIE for larger problems to demonstrate the excellent accuracy, iterative solver convergence and stability at interior resonance frequencies. 

\section{Integral Equations for Scattering Problems} %
\subsection{Single Field Formulations}
A typical scattering scenario is shown in Fig.~\ref{fig:scatterer}. 
A PEC scattering object is located inside the source domain $V_{\mathrm{i}}$, enclosed by the surface $S_0$ and surrounded by free space, i.e.\ the solution domain $V_{\mathrm{e}}$. 
The incident field is denoted for the case of an incident plane wave by 
\begin{alignat}{3}
\bm {E}^{\mathrm{inc}}(t,\bm {r})&=\,&&\bm {E}_{0} \,\mathrm {e} ^{\,\mathrm{j}(\omega t - \bm {k}_0 \!\cdot \bm {r})}\,,\\
\bm {H}^{\mathrm{inc}}(t,\bm {r})&=\,&&\frac{1}{\omega\mu_0}\bm {k}_0\times\bm {E}_0 \,\mathrm {e} ^{\,\mathrm{j}(\omega t - \bm {k}_0 \!\cdot \bm {r})}
\end{alignat}
with a, in the following suppressed, time dependency $\mathrm{e} ^{\,\mathrm{j}\omega t}$. 
The angular frequency is called $\omega$, $\bm {k}_0$ is the wavevector (as a consequence, the wavenumber is $\left|\bm {k}_0\right|=k_0$) and $\bm {r}$ is the observation point located on the surface~$S_0$.  
Furthermore, the normal vector on the surface~$\bm {n}(\bm {r})$ is defined to point into the solution domain~$V_{\mathrm {e}}$. 

The scattered fields~$\bm E^{\mathrm{sca}}(\bm {r})$ and~$\bm {H}^{\mathrm{sca}}(\bm {r})$ are a consequence of the presence of a scatterer. 
According to the Huygens principle, the radiation caused by the object can be described by equivalent sources on the surface~$S_0$ of the scatterer~\cite{Schelkunoff1936}, as depicted in Fig.~\ref{fig:scatterer}(b). 
However, the equivalent source representation is not unique. 
In particular, the field inside $V_{\mathrm {i}}$ varies with the different representations of the equivalent sources. 
A quite special choice of the equivalent sources are the so-called Love currents~\cite{Love1901}. 
They are related to the tangential field components on $S_0$. More specifically, the electric surface current densities~\cite{Schelkunoff1936,Love1901}
\begin{equation}
\bm J_{\mathrm {S}}(\bm r) = \bm n(\bm r) \times \bm H(\bm r)\label{eq:Love1}
\end{equation}
and the magnetic surface current densities 
\begin{equation}
\bm M_\mathrm {S} (\bm r)= \bm E (\bm r)\times \bm n(\bm r)\,. \label{eq:Love2}
\end{equation}
are defined via the tangential fields on the Huygens surface. 
To fulfill the boundary condition for PEC,
\begin{equation}
\bm n(\bm r) \times \bm E(\bm r) = \bm 0\,,\label{eq:PEC-Bound}
\end{equation} 
the Love currents are given by 
\begin{equation}
\bm J_\mathrm {S}(\bm r) = \bm n \times \bm H(\bm r)\,, \qquad \bm M_\mathrm S (\bm r)= \bm 0\,. \label{eq:Love-PEC}
\end{equation}

The most remarkable consequence of the Love currents is that the interior of the source domain becomes field-free, as shown in Fig.~\ref{fig:scatterer}(b). 
Then, the next step, see Fig.~\ref{fig:scatterer}(c), is to replace the PEC object with free space, where the zero field  inside $V_\mathrm i$ and the scattered fields in $V_\mathrm e$ remain unchanged.  
As a consequence, the Green's function of free space
\begin{equation}
G_0(\bm r,\bm r')=\frac{\mathrm{e}^{-\mathrm jk_0\left|\bm r -\bm r'\right|}}{4\uppi\left|\bm r -\bm r'\right|}
\end{equation}
is employed to determine the radiated fields of the equivalent sources, where $\bm r'$ represents the source coordinate. 

In general, we can conclude that the equivalent sources will model the scattered fields in a way that their superposition with the incident field~\cite{peterson1998computational}
\begin{alignat}{3}
\bm n(\bm r) \times \bm E^\mathrm{inc}(\bm r)\,+\,&\bm n(\bm r) \times \bm E^\mathrm{sca}(\bm r)&&=\bm n(\bm r) \times \bm E(\bm r)\,,\label{eq:EFIE}\\
\bm n(\bm r) \times \bm H^\mathrm{inc}(\bm r)\,+\,&\bm n(\bm r) \times \bm H^\mathrm{sca}(\bm r)&&=\bm n(\bm r) \times \bm H(\bm r)\label{eq:MFIE}
\end{alignat}
on $S_0$ fulfills the boundary condition~\eqref{eq:PEC-Bound} on the PEC object.
If the scattered field is expressed in terms of equivalent Love surface current densities, we obtain the EFIE~\eqref{eq:EFIE} as
\begin{equation}
\bm n \times \bm E ^\mathrm{inc} =\mathrm j k_0 Z_0 \mathcal T \bm J_\mathrm S\label{eq:EFIE-anal}
\end{equation}
and the MFIE~\eqref{eq:MFIE} as
\begin{equation}
\bm n \times \bm H ^\mathrm{inc} =\mathcal M \;\!\bm J_\mathrm S = \big[\mathcal I/2-\mathcal K \big]\bm J _ \mathrm S\,,\label{eq:MFIE-anal}
\end{equation}
where~$Z_0$ represents the free-space wave impedance. 
The equations are abbreviated by the use of the EFIE operator
\begin{multline}
\mathcal T \;\! \bm J _ \mathrm S \coloneq  \bm n(\bm r) \times 
\oiint\nolimits_{S_0}G_0(\bm r,\bm r') \bm J_\mathrm{S}(\bm r)\mathop{}\!\mathrm{d}s'
\\
+\bm n(\bm r) \times
\frac{1}{k_0^2}\bm\nabla\oiint\nolimits_{S_0}G_0(\bm r,\bm r')\bm\nabla'\cdot\bm J_\mathrm{S}(\bm r)\mathop{}\!\mathrm{d}s'
\end{multline}
and the MFIE operator $\mathcal M$ consisting of the identity operator
\begin{equation}
\mathcal I \;\! \bm J _ \mathrm S\coloneq  \bm J _ \mathrm S
\end{equation}
and the MFIE integral operator
\begin{equation}
\mathcal K\;\! \bm J _ \mathrm S\coloneq  \bm n(\bm r) \times 
\oiint\nolimits_{S_0}\bm\nabla G_0(\bm r,\bm r')\times\bm J_\mathrm{S}(\bm r')\mathop{}\!\mathrm{d}s'\,.
\end{equation}
Omitting the Love condition~\eqref{eq:Love-PEC} changes the EFIE to 
\begin{equation}
\bm n \times \bm E ^\mathrm{inc} =\mathrm j k_0 Z_0 \mathcal T \bm J_\mathrm{S} + \big[\mathcal I/2+\mathcal K \big]\bm M_\mathrm S\,,\label{eq:EFIE-generic}
\end{equation}
including both electric and non-vanishing magnetic surface current densities. 
This solution can produce the exactly same fields in $V_\mathrm e$, however, the equation is under-determined with twice as many unknowns as necessary. 
Furthermore, this version of the EFIE might give non-zero fields in $V_\mathrm i$ for the superposition of the incident and the scattered fields. 

\subsection{Combined Field and Source Integral Equations}
As mentioned in the introduction, both EFIE and MFIE suffer from interior resonances. 
Hence, a superposition of the two equations to the combined field integral equation (CFIE)~\cite{mautz1979h}
\begin{multline}
\alpha \mathrm j k_0 Z_0 \mathcal T \bm J_\mathrm{S} + (1-\alpha)Z_0\, \bm n \times \mathcal M \;\! \bm J_\mathrm S \\=\alpha\bm n \times \bm E^\mathrm{inc}+(1-\alpha)Z_0\bm n \times \bm n \times \bm H^\mathrm{inc}\label{eq:CFIE}
\end{multline}
with the combination factor $\alpha$ eliminates these resonances and leads to a unique solution. 
The common notation of a rotated MFIE is a convention for the choice of rotated testing functions as compared to the EFIE. 
In principle, EFIE and MFIE testing functions can, and should, be chosen separately for optimal performance as it is done in mixed discretization formulations. 

The CSIE forms an alternative to the CFIE with Love currents. 
In this formulation, the magnetic currents~\cite{Mautz1979}
\begin{equation}
\bm M_\mathrm{S,CS}=\beta Z_0\bm n \times \bm J_\mathrm{S,CS}\label{eq:CS-strong}
\end{equation}
are defined as point-wise rotated, scaled versions of the electric currents. 
These combined sources, also known as Huygens radiators, show a directive radiation characteristic into the solution domain and, thus, approximately suppress the scattered field inside the source region. 
With a weighting of $\beta=1$, the radiation contributions of electric and magnetic currents are balanced, i.e.\ this choice corresponds to a CFIE with $\alpha=0.5$.  

Based on the EFIE~\eqref{eq:EFIE-generic} for PEC scatterers, the CSIE~\cite{Mautz1979}
\begin{equation}
\mathrm j k_0 Z_0 \mathcal T \bm J_\mathrm{S} + \beta Z_0 \big[\mathcal I/2+\mathcal K \big]\bm n \times \bm J_\mathrm {S,CS} = \bm n \times \bm E^\mathrm{inc}\label{eq:CS-EFIE}
\end{equation}
 with strong form CS condition~\eqref{eq:CS-strong} is defined. 
 The solution of~\eqref{eq:CS-EFIE} is unique and interior resonances are avoided. 
Furthermore, it has convergence benefits as compared to the EFIE with Love currents: 
The combined sources minimize interactions between equivalent sources by suppressing the radiation through the scatterer. 
This is not the case for independent electric and magnetic currents which radiate equally to the outside and through the scatterer. 

\section{Discretization} %
\subsection{Love Current Formulations}
In the following, the integral equations will be discretized and evaluated on a triangular mesh. 
The electric and magnetic surface current densities are modeled by RWG basis functions~$\bm \beta(\bm r)$ on pairs of adjacent triangles according to
\begin{equation}
\bm {J}_{\mathrm {S}}(\bm {r})= \sum\nolimits_{n=1}^N \vecelem {i}_n \bm {\beta}_n(\bm {r})\,,\quad
\bm {M}_{\mathrm {S}}(\bm {r})= \sum\nolimits_{n=1}^N \vecelem {v}_n \bm {\beta}_n(\bm {r})\,,
\end{equation}
where $\vecelem {i}_n$ and $\vecelem {v}_n$ form the vectors $\vec {i}$ and $\vec {v}$ with $N$ unknown expansion coefficients each.

Based on the surface and source discretization, the incident and scattered fields are evaluated on the surface by weighting with appropriate testing functions. 
In the ideal case, the testing functions should be chosen from the dual space of the quantity to evaluate. 
For the  operator $\mathcal T$, it is well-known that $\bm \alpha =\bm n \times\bm\beta$ functions are the proper testing functions in the dual space of $\bm n \times \bm E$. 
Alternatively, $\bm\beta$ testing functions can be used for $\bm n \times \bm n \times \bm E$. 
This leads to the discretized version of the EFIE~\eqref{eq:EFIE-anal}
\begin{equation}
\mathrm j k_0 Z_0 \mat T_{\bm\beta} \vec i=\mathrm j k_0 Z_0 \Big[ \mat B_{\bm\beta} +  \frac{1}{k_0^2} \mat C_{\bm\beta}\Big]\vec i = \vec e_{\bm\beta}\label{eq:EFIE-discrete}
\end{equation}
for PEC objects and Love currents. $\mat B_{\bm\beta}$ and $\mat C_{\bm\beta}$ are matrices with the elements
\begin{equation}
 \big[\mat B_{\bm\beta}\big]_{mn}  =  \iint\nolimits_{S_m}\!\bm \beta_m(\bm r)\cdot\!\iint\nolimits_{S_n}G_0(\bm r , \bm r')\bm\beta_n(\bm r')\mathop{}\!\mathrm{d}s'\mathop{}\!\mathrm{d}s\,,
\end{equation}
\begin{equation}
 \big[\mat{C}_{\bm\beta}\big]_{mn}  = \iint\nolimits_{S_m}\!\bm \beta_m(\bm r)\cdot\bm\nabla\!\iint\nolimits_{S_n}G_0(\bm r , \bm r')\bm\nabla'_\mathrm{S}\cdot\bm\beta_n(\bm r')\mathop{}\!\mathrm{d}s'\mathop{}\!\mathrm{d}s
\end{equation}
and $\vec{e}_{\bm\beta}$ is a column vector with the entries
\begin{equation}
\big[\vec{e}_{\bm\beta}\big]_m = \iint\nolimits_{S_m}\!\bm\beta_m(\bm r) \cdot\bm E ^\mathrm{inc}(\bm r)\mathop{}\!\mathrm{d}s\,.
\end{equation}

In the case of the MFIE, the situation is not as clear as for the EFIE. 
In the form given in~\eqref{eq:MFIE-anal}, as indicated by the identity operator, the RWG basis functions are mapped from their original space, $\mathcal H_\mathrm{div}^{-1/2}$, to the same space again.  
Therefore, the common choice of RWG testing functions does not fulfill dual-space testing of the~$\mathcal I$ and~$\mathcal K$ operators. 
One possibility is to utilize other functions, most mentionable the BC functions, as they reside in the dual space of $\mathcal H_\mathrm{div}^{-1/2}$, which is $\mathcal H_\mathrm{curl}^{-1/2}$. 
If RWG functions are considered, the natural way to test in the proper space would be with rotated RWG functions $\bm \alpha$. 
This leads, however, to a degenerate matrix with zeros on the main diagonal and, hence, to poor iterative solver convergence. 

Despite the widespread concerns about the testing with RWG functions, we discretize the MFIE~\eqref{eq:MFIE-anal} in the classical manner
\begin{equation}
\mat M \vec i =-\Big[\frac{1}{2} \mat G_{\,\bm\beta\bm\beta}+ \mat K_{\,\bm\alpha}\Big]\vec i = \vec h_{\,\bm\alpha}\,.\label{eq:disc-MFIE}
\end{equation}
The matrices are computed element-wise as
\begin{equation}
 \big[\mat{G}_{\,\bm\beta\bm\beta}\big]_{mn}  =  \iint\nolimits_{S_m}\!\bm \beta_m(\bm r)\cdot\bm\beta_n(\bm r)\mathop{}\!\mathrm{d}s\,,\label{eq:Gbb}
\end{equation}
\begin{equation}
 \big[\mat{K}_{\bm\alpha}\big]_{mn}   = \iint\nolimits_{S_m}\!\bm \alpha_m(\bm r)\cdot\!\!\iint\nolimits_{S_n}\!\bm\nabla G_0(\bm r , \bm r')\times\bm\beta_n(\bm r')\mathop{}\!\mathrm{d}s'\mathop{}\!\mathrm{d}s
\end{equation}
and, similarly, the entries of the tested magnetic field vector
\begin{equation}
\big[\vec{h}_{\bm\alpha}\big]_m = \iint\nolimits_{S_m}\!\bm \alpha_m(\bm r) \cdot\bm H ^\mathrm{inc}(\bm r)\mathop{}\!\mathrm{d}s\,.
\end{equation}

The EFIE and the MFIE are then combined according to~\eqref{eq:CFIE} and the CFIE
\begin{equation}
\big[\alpha\mathrm j k_0Z_0\mat T_{\bm \beta}+(1-\alpha)Z_0\mat M\big]\vec i = \alpha \vec e_{\bm\beta} + (1-\alpha)Z_0 \vec h_{\bm\alpha}
\end{equation}
is obtained.

\subsection{Combined-Source Formulation}
\begin{figure*}[t]
\centering
 \begin{minipage}[b]{0.245\linewidth}
 \centering%
  \subfloat[]{%
 \includegraphics[]{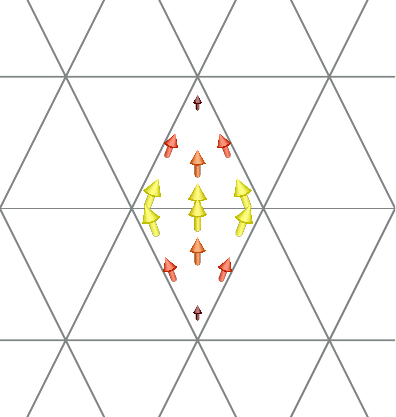}%
   }%
 \end{minipage}
 \hfill
 \begin{minipage}[b]{0.245\linewidth}
 \centering%
  \subfloat[]{%
 \includegraphics[]{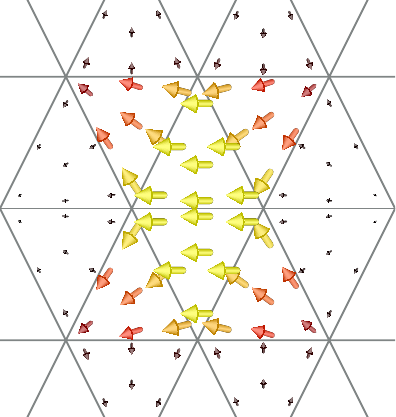}%
   }%
 \end{minipage}
 \hfill
 \begin{minipage}[b]{0.245\linewidth}
 \centering%
  \subfloat[]{%
 \includegraphics[]{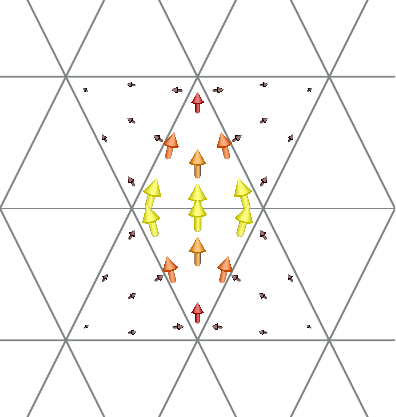}%
   }%
 \end{minipage}
 \hfill
 \begin{minipage}[b]{0.245\linewidth}
 \centering%
  \subfloat[]{%
 \includegraphics[]{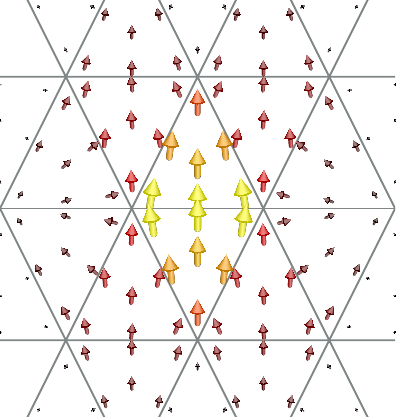}%
   }%
 \end{minipage}
 \caption{Surface current density discretization with RWG basis functions on a triangular mesh. (a) Standard RWG basis function. (b) Weak rotation $\mat{G}_{\bm \beta \bm \beta}^{-1}\mat{G}_{\bm \beta\bm \alpha}$ applied for this very RWG basis function. (c) Evaluated identity operator of the classical MFIE, i.e.\ one column of the Gram matrix $\mat G_{\bm \beta \bm \beta}$. (d) Evaluated identity operator of the new WMFIE-1, i.e.\ one column of the matrix $\mat R$. \label{fig:quiver}}
 \vspace*{-0.2cm}
\end{figure*}

Before going into an improved discretization of the identity operator within the MFIE, let us recall the CSIE with weak-form CS condition. 
This key idea will then be transferred to the MFIE. 
For the discretization of the EFIE~\eqref{eq:EFIE-generic}, RWG basis functions are chosen for both electric and magnetic surface current densities. Standard testing leads to the equation system
\begin{equation}
\Big[-\frac{1}{2} \mat G_{\bm\beta\bm\alpha} +  \mat K_{\bm\beta}\Big]\vec v +
\mathrm j k_0 Z_0 \mat T_{\bm\beta} \vec i = \vec e_{\bm\beta}\,,\label{eq:EFIE-CS-disc}
\end{equation}
 including the matrices $\mat{G}_{\bm\beta\bm\alpha}$ and $\mat{K}_{\bm\beta}$ with the entries
\begin{equation}
 \big[\mat{G}_{\bm\beta\bm\alpha}\big]_{mn}  =  \iint\nolimits_{S_m}\!\bm\beta_m(\bm r)\cdot\bm \alpha_n(\bm r)\mathop{}\!\mathrm{d}s\,,\label{eq:Gba}
\end{equation}
\begin{equation}
 \big[\mat{K}_{\bm\beta}\big]_{mn}   = \iint\nolimits_{S_m}\!\bm \beta_m(\bm r)\cdot\!\!\iint\nolimits_{S_n}\!\bm\nabla G_0(\bm r , \bm r')\times\bm\beta_n(\bm r')\mathop{}\!\mathrm{d}s'\mathop{}\!\mathrm{d}s\,.
\end{equation}
The important feature is the discretization of the CS condition~\eqref{eq:CS-strong}. Since the div-conforming~$\bm\beta$ functions cannot be used to directly implement the strong-form CS condition---only the curl-conforming~$\bm\alpha$ functions with hyper-singular line charges could do so---,~\eqref{eq:CS-strong} is transferred into a weak-form condition
\begin{equation}
\mat G_{\,\bm\beta\bm\beta}\vec v =\beta Z_0 \mat G_{\bm\beta\bm\alpha}\vec i\label{eq:CS-weak}
\end{equation}
by testing it with RWG functions~\cite{Kornprobst_2017}. 
This weak-form CS condition performs the magnetic current vector rotation by $90^\circ$ globally within the RWG basis. 
Considering only one element of the electric current vector, the result of the weak-form rotation is depicted in Fig.~\ref{fig:quiver}(b). 
It is observed that the orthogonal current flow is modeled well inside the original pair of triangles. 
In addition, the current is distributed over neighboring triangles.

\section{Accurate Identity Operator Discretization\\in the MFIE} 
One possible reason for the inaccuracy of the MFIE is the strong singularity of the identity operator~\cite{ErguelGuerel2009,Yan2011identity}.
The presented CSIE mitigates this problem by the discretization of the weak-form side condition~\eqref{eq:CS-weak} and, thus, shows a higher accuracy than the classical MFIE~\cite{Kornprobst_2017}. 
From an IE operator point of view, the strong form CS condition, as well as the identity operator, are highly singular operators. 
This becomes obvious by writing the Gram matrix entries~\eqref{eq:Gbb} and~\eqref{eq:Gba} as 
\begin{equation}
 \big[\mat{G}_{\bm\beta\bm\beta}\big]_{mn}  =  \iint\nolimits_{S_m}\!\bm\beta_m(\bm r)\cdot \iint\nolimits_{S_n}\!\bm \beta_n(\bm r')\,\updelta(\bm r,\bm r')\mathop{}\!\mathrm{d}s'\mathop{}\!\mathrm{d}s\,,
\end{equation}
\begin{equation}
 \big[\mat{G}_{\bm\beta\bm\alpha}\big]_{mn}  =  \iint\nolimits_{S_m}\!\bm\beta_m(\bm r)\cdot\iint\nolimits_{S_n}\!\bm \alpha_n(\bm r')\,\updelta(\bm r,\bm r')\mathop{}\!\mathrm{d}s'\mathop{}\!\mathrm{d}s\,,
\end{equation}
i.e. with a Dirac delta integral kernel. 
The high spectral content of these operators requires careful discretization to obtain an isotropic averaging low-pass effect~\cite{ErguelGuerel2009,Davis2004higherorder,Warnick2007regularization,Yan2011identity}. 
However, the integration of the weak-form CS condition~\eqref{eq:CS-weak} inside the discretization of the MFIE operator $\mathcal M$ seems not to be possible in a directly obvious way, since it approximates a $90^\circ$ vector rotation. 
Therefore, we introduce fictitious currents~$\bm{\hat{J}}_{\mathrm S}$ and~$\bm{\tilde{J}}_{\mathrm S}$ and perform the rotation as in~\eqref{eq:CS-strong} twice 
\begin{equation}
\bm{\hat{J}}_S
=-\bm n \times \bm{\tilde{J}}_{\mathrm S}=-\bm n \times\bm n \times \bm J_\mathrm{S}=\bm J_\mathrm{S}\,,
\end{equation}
where~$\bm{\hat{J}}_{\mathrm S}$ is equal to~$\bm J_\mathrm{S}$ and~$\bm{\tilde{J}}_{\mathrm S}$ is a $90^\circ$-rotated intermediate version of the current. 
Then, it is possible to discretize both strong-form rotations, as it is done in the CS condition discretization, to obtain the two systems of equations
\begin{equation}
\mat G_{\,\bm\beta\bm\beta}\bm\hat{\vec \i\hspace*{0.065cm}}\hspace*{-0.065cm}
=-\mat G_{\bm\beta\bm\alpha}\bm\tilde{\vec \i\hspace*{0.075cm}}\hspace*{-0.075cm}\,,\qquad
\mat G_{\,\bm\beta\bm\beta}\bm\tilde{\vec \i\hspace*{0.075cm}}\hspace*{-0.075cm}=\mat G_{\bm\beta\bm\alpha}\vec i
\,,
\end{equation}
where the discretized current~$\bm\hat{\vec \i\hspace*{0.065cm}}\hspace*{-0.065cm}$ is employed to approximate the current~$\vec i$. 
This becomes more obvious if the approximation step is summarized in the matrix
\begin{equation}
\mat R =-\mat{G}_{\bm\beta\bm\beta}^{-1}\mat G_{\bm\beta\bm\alpha}\mat{G}_{\bm\beta\bm\beta}^{-1}\mat G_{\bm\beta\bm\alpha}\,,
\end{equation}
which now contains two subsequent weak-form vector rotations by $90^\circ$ and a multiplication by $-1$.
Since the matrix $\mat G_{\bm\beta\bm\alpha}$ employed for the weak-form rotation~\eqref{eq:CS-weak} is singular, the CSIE cannot be employed for large weighting factors of the magnetic currents. 
Similarly, $\mat R$ inherits the null-space of the matrix $\mat G_{\bm\beta\bm\alpha}$ and, hence, cannot be employed solely as weak-form identity operator discretization. 
In order to avoid this very problem in the MFIE, we  combine $\mat R$ with the strong-form identity matrix~$\mat I$ to a weak-form basis functions transformation matrix
\begin{equation}
\mat W_\gamma = \gamma \mat {I} + (1-\gamma)\mat R\,,
\end{equation}
where $\gamma$ is a  weighting factor with $0\leq\gamma\leq1$. 

There are several possibilities to utilize such a weak-form relation within the discretized MFIE~\eqref{eq:disc-MFIE} and several new MFIE formulations with weak-form rotation (WMFIE) are derived. 
The first possibility concentrates on the discretized identity operator. 
This is motivated by the knowledge that the identity operator is the main cause of the accuracy problems. 
As opposed to the CSIE, the $\mathcal K$ operator does not seem to need a special treatment. 
Nevertheless, this is also investigated in the following. 
The WMFIE-1 formulation is given as
\begin{equation}
-\Big[\frac{1}{2} \mat G_{\,\bm\beta\bm\beta}\mat W_\gamma + \mat K_{\,\bm\alpha}\Big]\vec i = \vec h_{\bm\alpha}\,,\label{eq:wmfie1}
\end{equation}
 applying~$\mat W_\gamma$ to the basis functions of the Gram matrix $\mat G_{\bm\beta\bm\beta}$. 

Alternatively, it is possible to multiply the matrix $\mat W_\gamma$ to the whole MFIE matrix. 
Then, the MFIE formulation is more similar to the one of the CSIE, since the $\mathcal K$ operator is also involved. 
This investigation can be seen as an empirical study of the effect of the weak-form identity operator. 
Multiplying $\mat W_\gamma$ from the right gives the WMFIE-2 formulation
\begin{equation}
\mat M \mat W_\gamma\vec i = \vec h_{\bm\alpha}
\end{equation}
and implementing a weak-form testing functions transformation matrix (from the left) results in the WMFIE-3 formulation
\begin{equation}
\big[\gamma \mat {I} + (\gamma-1)\mat G_{\bm\beta\bm\alpha}\mat{G}_{\bm\beta\bm\beta}^{-1}\mat G_{\bm\beta\bm\alpha}\mat{G}_{\bm\beta\bm\beta}^{-1}\big]\mat M \vec i = \vec h_{\bm\alpha}\,.
\end{equation}
\begin{equation}
\mat{G}_{\bm\beta\bm\beta}\mat W_\gamma\mat{G}_{\bm\beta\bm\beta}^{-1}\mat M \vec i = \vec h_{\bm\alpha}\,.
\end{equation}
The weak-form testing functions transformation matrix employed solely for $\mat{G}_{\bm\beta\bm\beta}$ results again in WMFIE-1~\eqref{eq:wmfie1}.

To gain further insight,  the basis transformation scheme of the WMFIE-1 formulation and its representation of the discretized identity operator in the two extreme cases $\gamma=0$ and $\gamma=1$ is investigated in Fig.~\ref{fig:quiver}(c) and (d). 
The important property of the identity operator representation with $\gamma =0$, see Fig.~\ref{fig:quiver}(d), is that the original basis function is distributed over several triangles. 
This is a quite similar effect as the testing with BC functions, which also enlarges the testing area to more than the source triangle~\cite{Yan2011identity}.

\begin{figure}[pt]
 \begin{minipage}[b]{0.49\linewidth}
 \centering%
  \subfloat[]{%
 \includegraphics[]{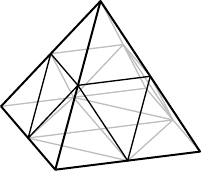}%
   }%
 \end{minipage}%
 \hfill
 \begin{minipage}[b]{0.49\linewidth}
 \centering%
  \subfloat[]{%
 \includegraphics[]{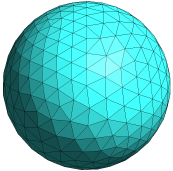}%
   }%
 \end{minipage}%
 \caption{Scatterers for analyzing the weighting factor $\gamma$. (a) A small pyramid. (b)~A~sphere analyzed below its first resonance. \label{fig:scatt1}}
 \vspace*{-0.2cm}
\end{figure}

The influence of the weighting parameter $\gamma$ has to be analyzed in the four presented versions of the WMFIE. 
Two different scattering objects, depicted in Fig.~\ref{fig:scatt1}, are considered: a small sharp-edged pyramid as an example of a sharp-edged object with coarse discretization and a small sphere.

\begin{figure}[pt]
 \begin{minipage}[b]{\linewidth}
 \centering%
  \subfloat[]{%
 \includegraphics[]{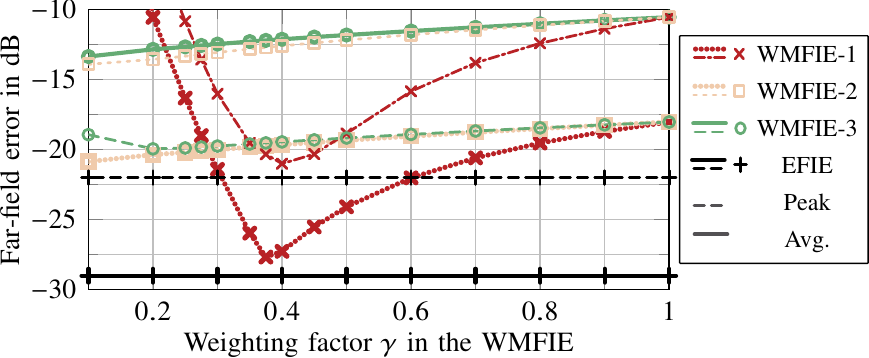}%
   }%
 \end{minipage}%
 \\[1.5ex]
 \begin{minipage}[b]{0.49\linewidth}
 \centering%
  \subfloat[]{%
 \includegraphics[]{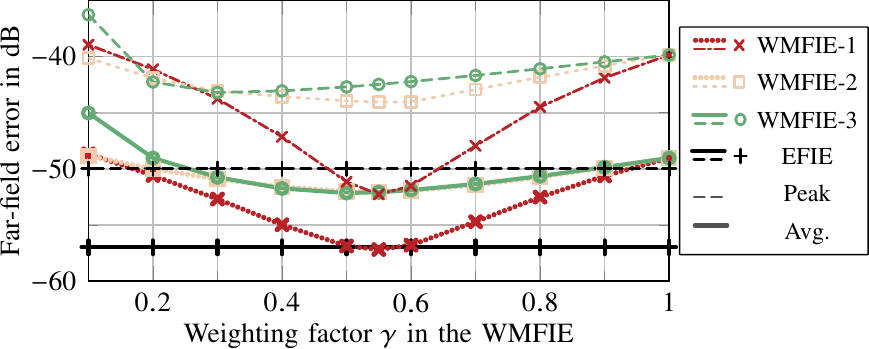}%
   }%
 \end{minipage}%
 \caption{Analysis of the weighting factor $\gamma$ for the formulations WMFIE-1 to WMFIE-3. (a) Results for the small pyramid. (b) Results for a sphere with 1\,m diameter. \label{fig:gamm}}
 \vspace*{-0.11cm}
\end{figure}

Reference solutions are the EFIE on a refined mesh for the pyramid and Mie series expansion for the sphere. 
The far-field accuracy of the bistatic radar cross section (RCS) over~$\gamma$ is shown in Fig.~\ref{fig:gamm}(a) for the pyramid, which is discretized with~$\lambda/10$ mean edge length. 
For $\gamma\rightarrow 1$, the standard MFIE (and the corresponding error level) is obtained in all three formulations. 
For the MFIE, an error of $-18$\,dB and a maximum error of $-11$\,dB are observed. 
The EFIE, with a very accurate solution on the same mesh, shows a maximum error of $-22$\,dB and an average error of  $-29$\,dB.  
Looking at the results of the WMFIE for small~$\gamma$ values, the null space of~$\mat G_{\bm\beta\bm\alpha}$ deteriorates the solution quality. 
This seems to be more severe for WMFIE-1. 
Nevertheless, this version shows the best error at around $\gamma=0.4$, almost achieving the EFIE accuracy. 
The accuracies of WMFIE-2 and WMFIE-3, which apply the weak-form rotation to the whole MFIE matrix~$\mat M$, are not as good as for WMFIE-1, but still better than the classical MFIE. 
This confirms our initial assumption that the accuracy problem of the MFIE originates from the discretized identity operator. 
The weak-form discretization of the $\mathcal K$ operator is neither meaningful nor does it give an improvement in performance as large as the weak-form representation of solely the $\mathcal I$ operator.%

The second example is a sphere with 1\,m diameter and a frequency of 200\,MHz, i.e.\ below the first interior resonance. 
Simulations are performed  with a rather fine discretization of $\lambda/14$ mean triangle edge length. The error for different weighting factors is shown in Fig.~\ref{fig:gamm}(b). 
Again, the WMFIE solution is severely deteriorated for $\gamma \rightarrow 0$ by the null space of $\mat G_{\bm\beta\bm\alpha}$, and for  $\gamma \rightarrow 1$ the standard MFIE accuracy is retained, i.e.\ a maximum error of $-47$\,dB and an average error of  $-40$\,dB.  
For comparison, the EFIE maximum and average errors are $-50$\,dB and $-57$\,dB, respectively. 
The \mbox{WMFIE-1} achieves quite similar accuracy levels with $\gamma\approx0.5$. 
The other formulations show worse accuracy, especially regarding the maximum error, as it was already the case for the first example. 
Again, the new weak-form discretization of the $\mathcal I$ operator improves the accuracy much more than
employing it for the $\mathcal M$ operator.

We conclude that the formulations with $\mat W_\gamma$ applied only to the Gram matrix of the RWG functions show a better  accuracy. 
In particular, the best formulation WMFIE-1 (in the following just WMFIE)
\begin{equation}
-\Big[\frac{1}{4} \mat G_{\,\bm\beta\bm\beta}-\frac{1}{4}\mat G_{\bm\beta\bm\alpha}\mat{G}_{\bm\beta\bm\beta}^{-1}\mat G_{\bm\beta\bm\alpha} + \mat K_{\,\bm\alpha}\Big]\vec i = \vec h_{\bm\alpha}\label{eq:WMFIE-final}
\end{equation} 
is also the most efficient one to implement, since it only requires one inversion of $\mat G _{\bm\beta\bm\beta}$, and not two as the other formulations. 
Unless otherwise stated, $\gamma=0.5$ is assumed. 
For this WMFIE, only the identity operator representation is modified. 
In consequence, the integral operators of EFIE and MFIE can be efficiently stored in one CFIE matrix. 
Matrix compression techniques such as the multilevel fast multi-pole method (MLFMM) or the adaptive cross approximation are immediately compatible with the presented scheme and the extra memory for storing the Gram matrices~$\mat G_{\bm\beta\bm\alpha}$ and~$\mat G_{\bm\beta\bm\beta}$ grows only with linear complexity and is negligible. 
The solution time is almost unaffected, as the iterative-solver inversion of~$\mat G_{\bm\beta\bm\beta}$  
for each matrix vector product is very fast. 

\section{Analysis of the WMFIE} %
\subsection{Convergence to the Correct Solution for Mesh Refinement}
A very important property to examine is the convergence to the correct solution when refining the mesh. 
Therefore, two examples are studied at a frequency of 150\,MHz with plane wave illumination. 
The considered integral equations are the classically discretized EFIE~\eqref{eq:EFIE-discrete} and the MFIE~\eqref{eq:disc-MFIE}, as well as the CSIE as the combination of~\eqref{eq:EFIE-CS-disc} and~\eqref{eq:CS-weak}
with $\beta=10$ and the WMFIE~\eqref{eq:WMFIE-final}. 
A pure magnetic currents CSIE is not possible due to the singular matrix~$\mat G_{\bm\beta\bm\alpha}$ appearing twice in~\eqref{eq:EFIE-CS-disc} and~\eqref{eq:CS-weak}. 
 Thus, the large $\beta$ coefficient mimics a CFIE with larger weighting of the MFIE part, i.e.\ $\alpha=1/11$, and gives insight about the behavior of the $\mathcal M$ operator discretization of the CSIE.

\begin{figure}[!t]
 \centering%
  \subfloat[]{%
 \centering
 \includegraphics{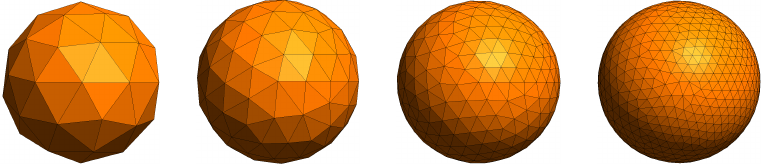}%
}\\\subfloat[]{%
 \includegraphics{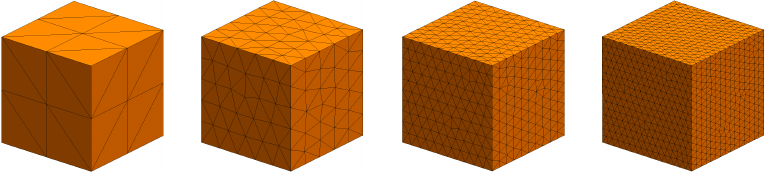}%
}\caption{Some of the models for the mesh refinement analysis, from coarse to fine mesh. (a)~Spheres. (b)~Square cubes.\label{fig:models2}}
\vspace*{-0.2cm}
\end{figure}
\begin{figure}[!t]
 \centering%
  \subfloat[]{%
 \includegraphics[]{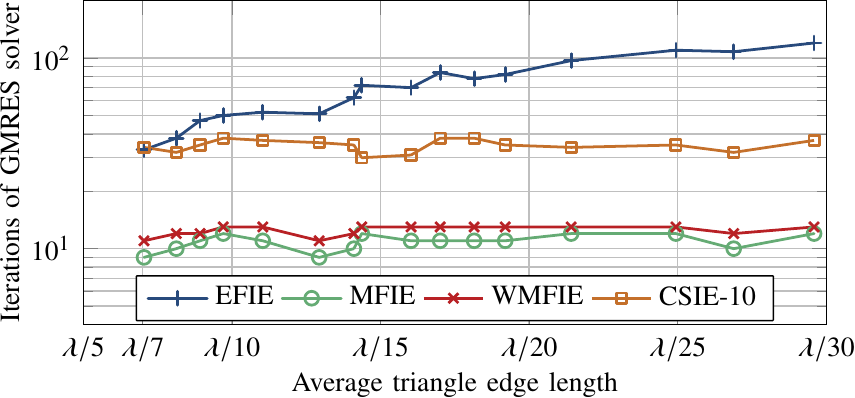}%
   }
   \\[1.5ex] 
  \subfloat[]{%
 \includegraphics[]{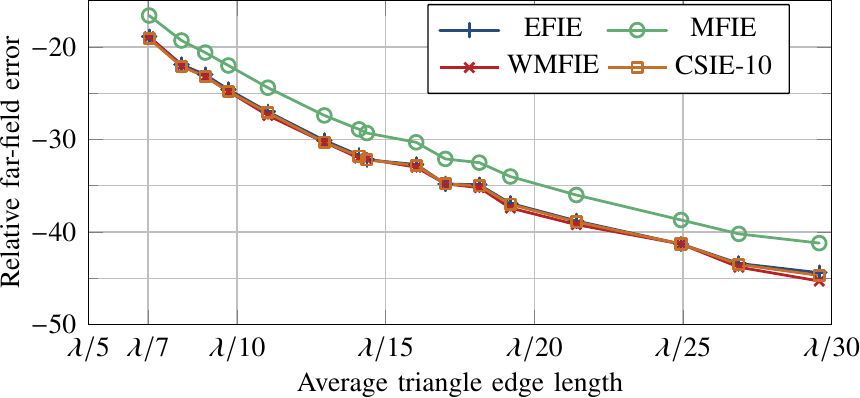}}%
   \\[1.5ex] 
  \subfloat[]{%
 \includegraphics[]{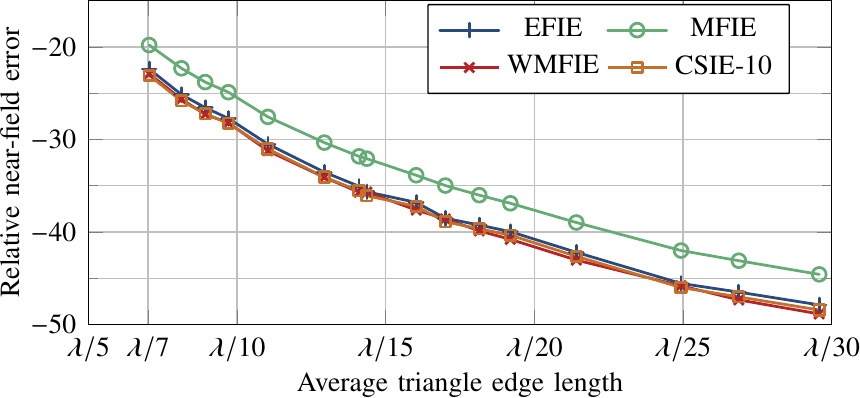}}%
 \caption{Mesh refinement analysis for a sphere with 1\,m diameter. (a)~Iterative solver convergence. (b)~Maximum relative error of the scattered far-field. (c)~Average relative near-field error on a 1.2\,m surrounding sphere.\label{fig:s-msh}}
\vspace*{-0.635cm}
\end{figure}
\begin{figure}[tp]
 \centering%
  \subfloat[]{%
 \includegraphics[]{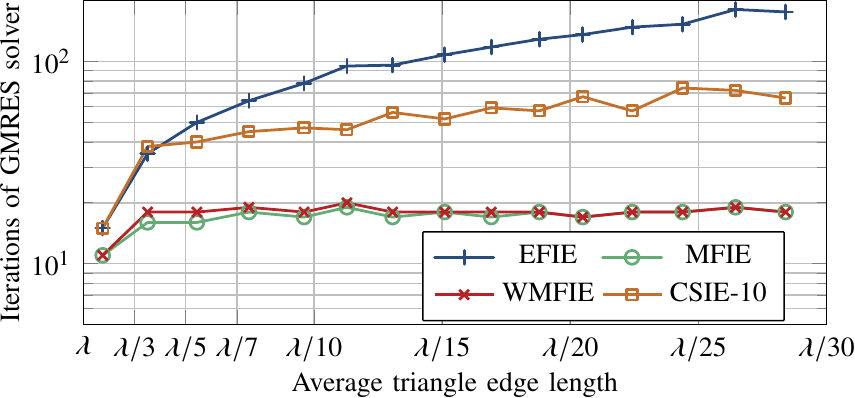}%
   }
   \\[0.9ex] 
  \subfloat[]{%
 \includegraphics[]{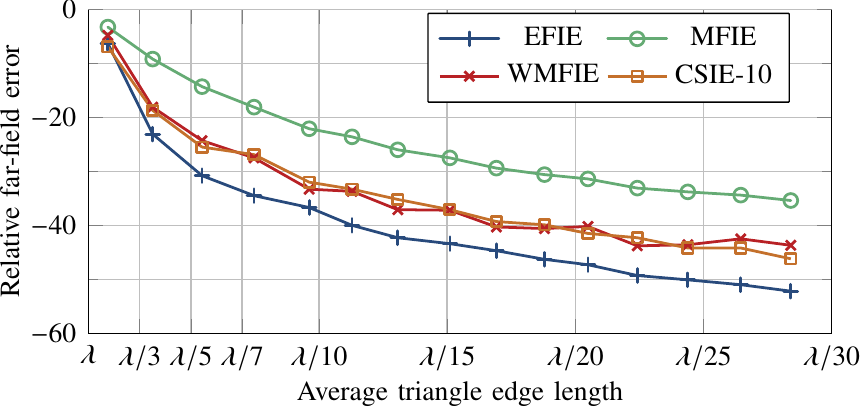}%
   }%
   \\[0.9ex] 
  \subfloat[]{%
 \includegraphics[]{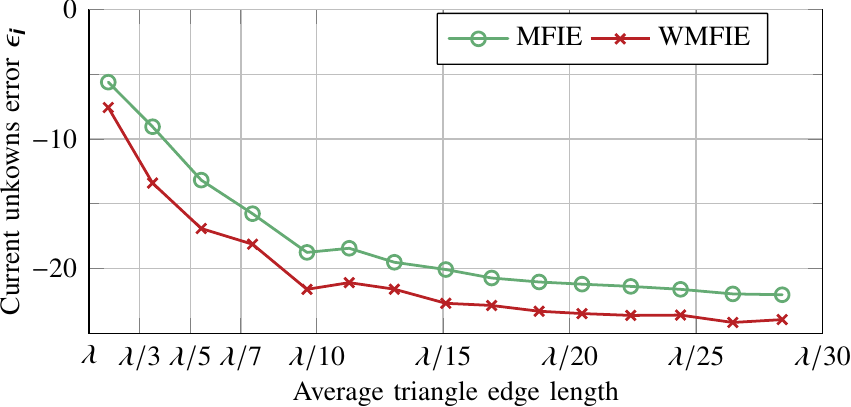}}%
 \caption{Mesh refinement analysis for a square cube. (a)~Iterative solver convergence. (b)~Maximum relative error of the scattered far-field. (c)~Comparison of MFIE and WMFIE solution vectors to the EFIE solutions.  \label{fig:c-mesh}}
\vspace*{-0.415cm}
\end{figure}

The first example is a sphere with a diameter of~1\,m, as shown in Fig~\ref{fig:models2}(a). The mesh refinement will be analyzed based on the average triangle edge length, with different simulation models for each data point. 
In Fig.~\ref{fig:s-msh}(a), the iterative solver convergence of a generalized minimum residual (GMRES) solver to a residual error of $10^{-4}$ is given. 
As expected, the EFIE shows a worse iteration behavior with refined mesh, whereas the MFIE converges with an approximately constant number of iterations. 
The WMFIE shows the same, fast convergence as the classical MFIE. The CSIE-10 shows a stable iterative solver convergence, however, worse than the MFIE. This has two reasons: First, the Gram matrix $\mat G_{\bm\beta\bm\alpha}$ has a null space and the MFIE part of the equation is singular. Second, to avoid the null space, the poorly converging EFIE is present with a weighting of $1/11$. 
The accuracy of the obtained solution is plotted in Fig.~\ref{fig:s-msh}(b) with Mie-series expansion as a reference. 
We computed the relative error 
\begin{equation}
\epsilon_{\vartheta,\varphi}(\vartheta,\varphi) = \left|\frac{ E^\mathrm{sca}_{\vartheta,\varphi,\mathrm{ref}}(\vartheta,\varphi)- E^\mathrm{sca}_{\vartheta,\varphi}(\vartheta,\varphi)}{\max\limits_{\vartheta, \varphi} \left(E_{\vartheta,\mathrm{ref}}(\vartheta,\varphi),E_{\varphi,\mathrm{ref}}(\vartheta,\varphi)\right)}\right|
\end{equation}
of the scattered far-field in several cut planes for the two polarizations~$E_\vartheta$ and~$E_\varphi$, where~$\bm E^\mathrm{sca}_\mathrm{ref}$ is the reference. 
Then, the maximum error~$\epsilon_\mathrm{max}=\max_{\vartheta, \varphi}\big(\epsilon_{\varphi}(\vartheta,\varphi),\epsilon_{\varphi}(\vartheta,\varphi)\big)$  is plotted as a measure of the accuracy. 
Other error measures like the average or the $l^2$-norm show the same behavior with somewhat other error values. 
 It is observed that the MFIE is about 3\,dB worse that the other IEs. 
 The same is observed for the average near-field error on a $1.2$\,m diameter surrounding sphere, see Fig.~\ref{fig:s-msh}(c). 
 The reason why the same accuracy is observed for CSIE-10, EFIE and WMFIE lies in the discretization of the sphere. 
There is a non-negligible discretization error of the curved surface and, furthermore, the discretized spheres are effectively  smaller than the ideal sphere considered in the Mie series expansion. This is obviously visible in Fig.~\ref{fig:models2}(a). 

One additional data point with a discretization of~$\lambda/167$ has been computed: Here, the EFIE converges within 250 iterations, while the WMFIE and MFIE still need about the same number of 12 and 11 iterations for convergence. 
More interestingly, the EFIE maximum error is at~$-69.9$\,dB, the WMFIE error at~$-65.2$\,dB and the MFIE~$-59.8$\,dB, 
i.e. the WMFIE is considerably better than the MFIE, while the EFIE shows here an even lower error. 
These error levels are achieved in comparison to the Mie series solution of a $0.9997$\,m sphere to correct for the smaller diameter of the discretized sphere.

The second exemplary scatterer, a square cube with 1\,m side length, is analyzed in the same way as the sphere. 
Lacking an analytical solution, the reference solution is taken from numerical simulation, i.e.\ the EFIE solution on the finest mesh with 3rd order expansion functions~\cite{Ismat2009MoM}, since it is expected that higher-order functions show a higher-order convergence to the correct solution with the mesh refinement. 
The iterative solver convergence to a residual error of $10^{-4}$ shown in Fig.~\ref{fig:c-mesh}(a) is similar to the first example: MFIE and WMFIE are equally well-conditioned and are not affected by the mesh refinement.
In the CSIE, the EFIE influence leads to a, though only slight, increase of the iteration count over the considered meshes. 
The accuracy analysis of the cube is more meaningful than for the sphere in a sense that the shape of the object is met exactly for any mesh.
Regarding the accuracy, the WMFIE and CSIE show about equal performance  and they are more  than 10\,dB better than the MFIE. However, the EFIE shows an even $5$\,dB to $10$\,dB better accuracy. A similar behavior has been observed for the mixed-discretization MFIE with BC testing functions~\cite{Cools_2011}. 
Furthermore, the relative error of solution vectors $\vec i_\mathrm{(W)MFIE}$ for the electric current unknowns are compared to the EFIE solution in Fig.~\ref{fig:c-mesh}(c), measured in the  $l^2$-norm according to 
\begin{equation}
\epsilon_{\vec i}=\frac{\lVert\vec i_\mathrm{EFIE}-\vec i_\mathrm{(W)MFIE}\rVert^2_2}{\lVert\vec i_\mathrm{EFIE}\rVert^2_2}\,.\label{eq:currenterror}
\end{equation}
During the mesh refinement, the WMFIE solution is always more accurate than the classically discretized MFIE solution.

\begin{figure}[t]
 \raggedleft%
  \subfloat[]{%
 \includegraphics[]{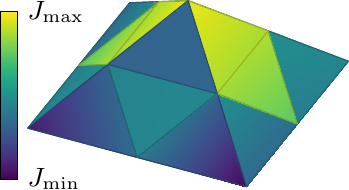}%
   }\hspace{0.75cm}%
  \subfloat[]{%
 \includegraphics[]{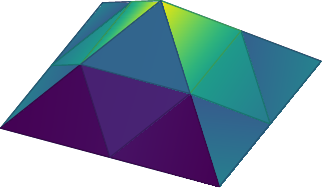}%
   }\hspace*{0.5cm}
   \\[1.5ex] 
  \subfloat[]{%
 \includegraphics[]{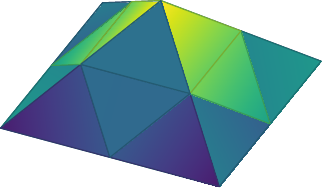}%
   }\hspace{0.75cm}%
  \subfloat[]{%
 \includegraphics[]{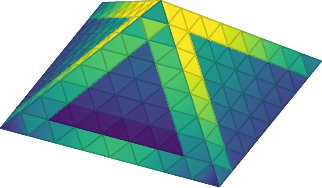}%
   }\hspace*{0.5cm}
 \caption{Absolute value of the surface current densities for different solutions obtained by: (a) EFIE, (b) MFIE, (c) WMFIE, (d) EFIE on refined mesh.  \label{fig:currents}}
\end{figure}%

\begin{figure}[tp]
 \centering%
 \includegraphics[]{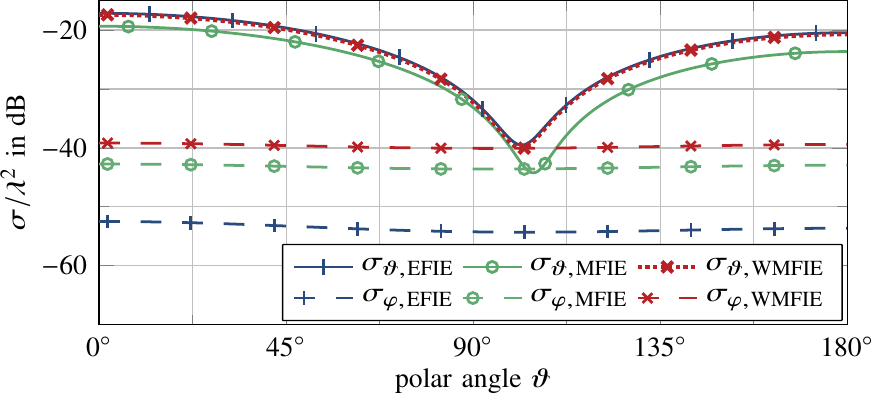}%
 \caption{Bistatic RCS of the small pyramid for the EFIE, MFIE and WMFIE solutions. The inaccuracy of the MFIE is clearly visible in the $\vartheta$ polarization.\label{fig:rcs-wmfie}}
\vspace*{-0.25cm}
\end{figure}

\subsection{Simulation Accuracy for Different Scattering Objects}

In the following, we analyze some simulation results of the WMFIE in more detail, in particular 
the near-field accuracy, and in the limiting case the surface current densities. 
This is examined for the example of the already discussed small pyramid. 
In Fig.~\ref{fig:currents}, the absolute value of the surface current density is plotted on the surface of the pyramid. The WMFIE and EFIE solutions give very similar results, whereas the MFIE solution differs quite a lot. 
Of course, the EFIE solution on a finer mesh looks quite different. 
Nevertheless, it is clear that the currents of the EFIE and WMFIE solutions are  more accurate than the MFIE on the coarse mesh. 

For the same example, the bistatic RCS for plane wave incidence is given in Fig~\ref{fig:rcs-wmfie}. 
The EFIE and WMFIE are quite close and show a maximum error as compared to the EFIE solution on a refined mesh of $-21.7$\,dB and $-18.8$\,dB, respectively. 
The MFIE exhibits a much larger error of $-10.6$\,dB. 
It is important to note that the far-field error is dominated by the co-polarization,~$\sigma_\vartheta$ in Fig.~\ref{fig:rcs-wmfie}. 
A discrepancy is seen in the cross-polarization; however, its absolute level is so low that it is not relevant in terms of accuracy. 

\begin{figure}[!tp]
 \centering%
  \subfloat[]{%
 \includegraphics[]{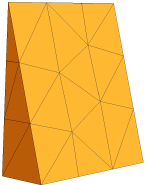}%
   }\hspace{1.5cm}%
  \subfloat[]{%
 \includegraphics[]{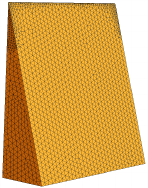}%
   }%
 \caption{Simulation models of the sharp wedge. (a) $\lambda/11$ mesh. \mbox{(b)~Refined~$\lambda/128$~mesh}. \label{fig:sharp_models}}
\end{figure}

\begin{figure}[!tp]
 \centering%
 \includegraphics[]{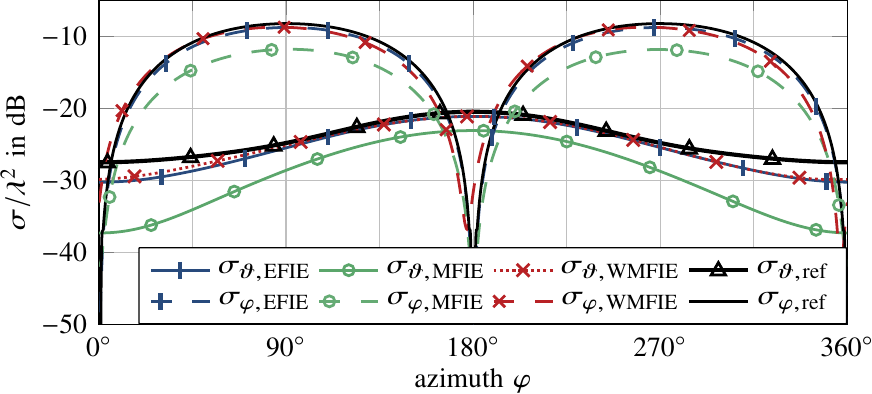}%
 \caption{Bistatic RCS of the sharp wedge. \label{fig:sharp}}
\vspace*{-0.25cm}
\end{figure}

\begin{figure}[!tp]
 \centering%
  \subfloat[]{%
 \includegraphics[]{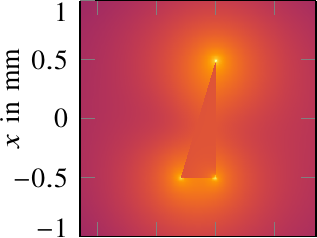}%
   }\hfill%
  \subfloat[]{%
 \includegraphics[]{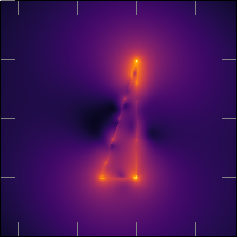}%
   }\hfill%
  \subfloat[]{%
 \includegraphics[]{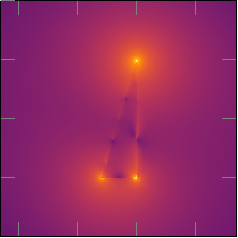}%
   }%
   \\
 \centering%
  \subfloat[]{%
 \includegraphics[]{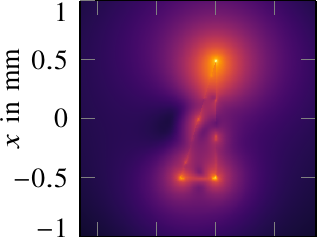}%
   }\hfill%
  \subfloat[]{%
 \includegraphics[]{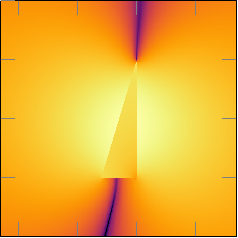}%
   }\hfill%
  \subfloat[]{%
 \includegraphics[]{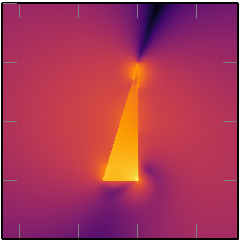}%
   }%
   \\
 \centering%
  \subfloat[]{%
 \includegraphics[]{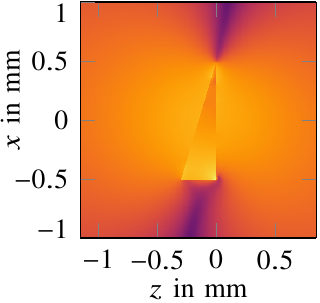}%
   }\hfill
  \subfloat[]{%
 \includegraphics[]{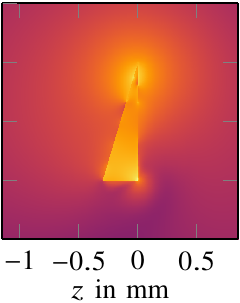}%
   }\hfill\hspace*{0.37cm}%
  \subfloat[]{%
 \includegraphics[]{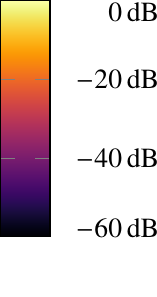}%
   }\hspace*{0.5cm}\mbox{}
 \caption{Scattered near field around and inside the sharp wedge: Normalized absolute field values and error vector magnitudes of the normalized fields. (a)~$E$ field reference solution. (b)~EFIE error. (c)~MFIE error. (d)~WMFIE error.  (e)~$H$ field reference solution. (f)~EFIE error. (g)~MFIE error. (h)~WMFIE error. (i)~Colorbar for all plots. \label{fig:sharp_nearfield}}
 \vspace*{-0.2cm}
\end{figure}

As a more extreme scenario, a sharp edge with 84 electric current unknowns has been simulated for plane-wave incidence with a mean edge length of $\lambda/11$. 
The EFIE solution on a refined mesh with $\lambda/128$ discretization is employed as a reference. 
Both models are depicted in Fig.~\ref{fig:sharp_models}. 
The RCS results are given in Fig.~\ref{fig:sharp}. 
The EFIE and WMFIE show maximum errors of about $-22$\,dB, while the MFIE is very inaccurate with a maximum error of $-9$\,dB. 
The scattered near field is analyzed in Fig.~\ref{fig:sharp_nearfield}. 
The absolute values of the electric and magnetic fields of the reference solution, as well as the spatial error distribution of EFIE, MFIE and WMFIE are depicted. 
It is observed that the WMFIE offers improvements for both electric and magnetic near fields, while it cannot achieve EFIE error levels completely. 
A word about the iterative solver convergence: The Gram matrix $\mat G_{\bm\beta\bm\beta}$ is somewhat ill-conditioned for such a mesh, it has a condition number of 810 (460 with diagonal preconditioning). The employed diagonally preconditioned conjugate gradient (CG) solver converges within 10 iterations to $10^{-10}$ for the weak-form rotation in the WMFIE, while the unpreconditioned GMRES for EFIE, WMFIE and MFIE converges within 63, 31, and 26 iterations to a residual of  $10^{-6}$. 
The computational cost for the weak-form rotation is overall more or less negligible, even for such a small problem.%


\subsection{Low-Frequency Behavior of the New WMFIE}
Another important aspect worth to be investigated is the low-frequency behavior of the new MFIE formulation. 
The classical RWG-tested MFIE shows a stable and good conditioning behavior for low frequencies, but the real part of the divergence of the current has erroneously a constant limit~\cite{zhang2003magnetic}. 
The reason is the testing of the~$\mathcal K$ operator~\cite{Bogaert_lowfrequencymfie_2011,Bogaert_lfmfie2_2014}, which is not changed in the presented formulation. 
Therefore, similar errors as for the classically discretized MFIE are expected for the presented formulation.

\begin{figure}[!tp]
 \centering%
  \subfloat[]{%
 \includegraphics[]{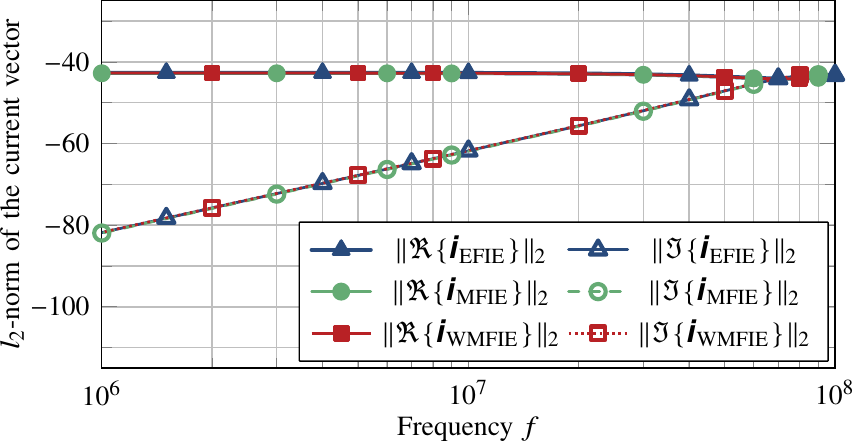}%
   }
   \\[1.5ex] 
  \subfloat[]{%
 \includegraphics[]{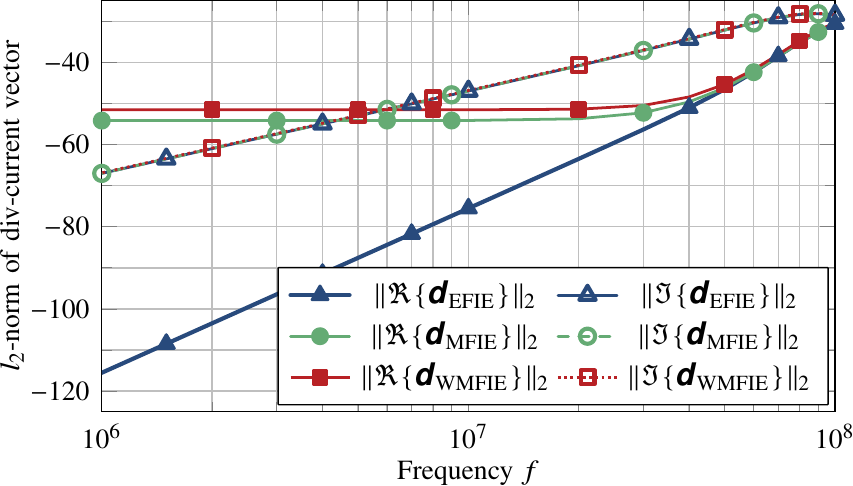}%
   }%
 \caption{Low-frequency behavior of classical and new MFIE formulations as compared to the correct EFIE solutions. (a)~Current unknowns behavior: everything seems to be fine. (b)~Divergence of the current: the constant but incorrect limit of the star/tree current is visible in the divergence.  \label{fig:lf}}
\end{figure}%

A 1\,m diameter PEC sphere with 126 unknowns has been simulated for decreasing frequencies. 
We investigate the real and imaginary parts of the current unkowns vector $\vec i$---containing RWG coefficients---as well as real and imaginary parts of the vector $\vec d$ containing the $\bm \nabla \cdot\bm J_\mathrm{S}$ coefficients for each triangle, see Fig.~\ref{fig:lf}. 
The divergence of the RWG current is associated with the charge inside the $k$th triangle%
\begin{equation}
\varrho_k = \frac{\mathrm j}{\omega} \bm \nabla \cdot\bm J_\mathrm{S},
\end{equation}
where $\bm \nabla \cdot\bm J_\mathrm{S}$ is  a constant times the coefficients~\cite{Rao_1982}. 
Furthermore, the divergence of the current has the same meaning as the star/tree current. 
As expected, a constant limit is observed in the real part of the divergence for all MFIE variations. 
Therefore, we conclude that the new MFIE is not suited for retrieving the star/tree current of low-frequency scattering problems.
The additional problems of not simply connected geometries 
are, naturally, also not tackled and not further discussed~\cite{ofluoglu2015magnetic,Cools_nullspace_2009}.

\begin{figure}[!tp]
 \centering%
 \includegraphics[]{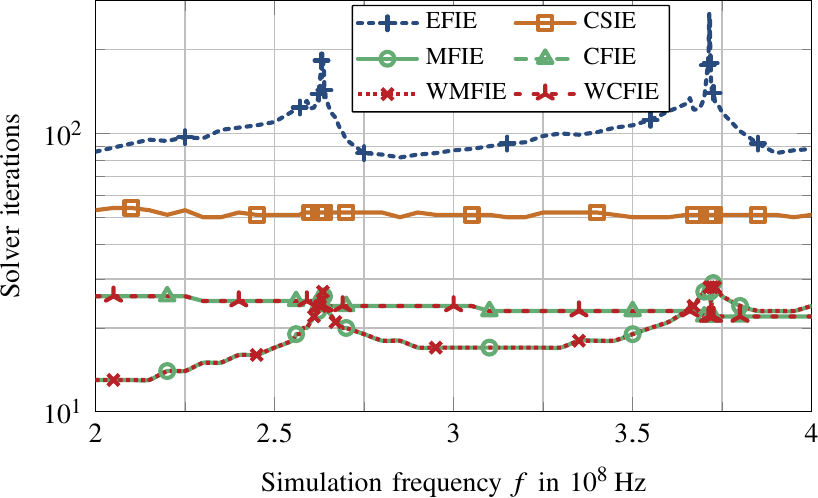}%
 \caption{GMRES iterative solver convergence at interior resonances for a 1\,m diameter sphere. \label{fig:s2-it}}
 \vspace*{-0.2cm}
\end{figure}

\begin{figure*}[!tp]
 \centering%
 \includegraphics[]{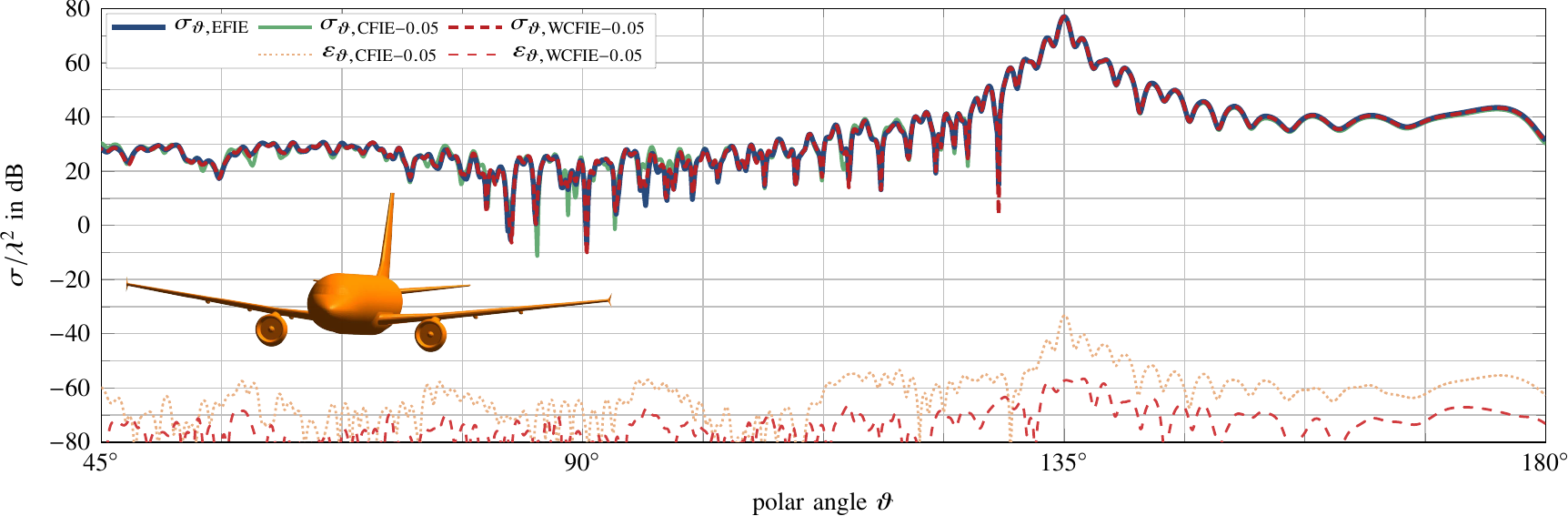}%
 \caption{Bistatic RCS of an A320 plane in the $\varphi=90^\circ$ cut plane. Comparison of EFIE, CFIE-0.05 and WCFIE-0.05 solutions.\label{fig:ab}}
\end{figure*}
\begin{figure}[!tp]
 \centering%
 \includegraphics[]{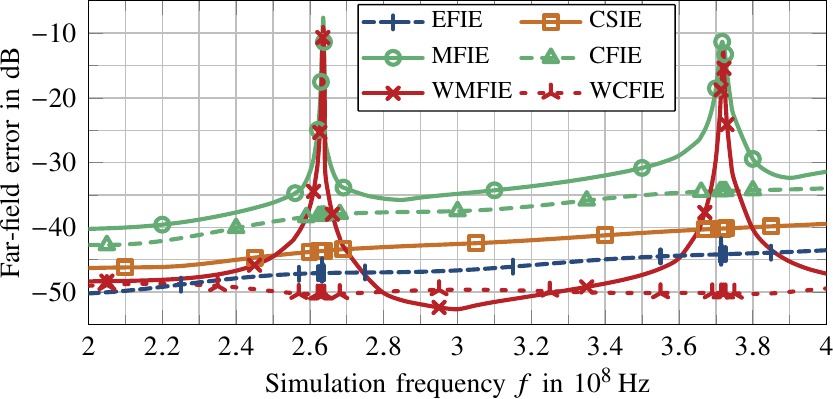}%
 \caption{Solution accuracy at interior resonances for a $1$\,m diameter sphere. \label{fig:s2-acc}}
\end{figure}
\begin{figure}[!tp]
 \centering%
 \includegraphics[]{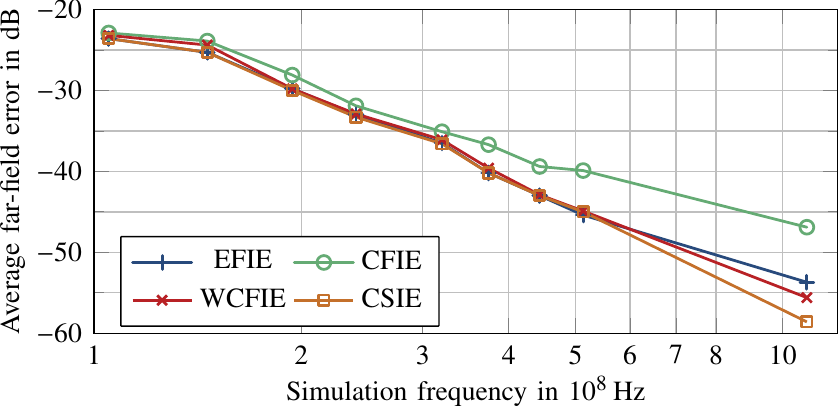}%
 \caption{Solution accuracy of a $1$\,m diameter sphere at varying frequencies. \label{fig:s-acc-freq}}
 \vspace*{-0.3cm}
\end{figure}

\section{Analysis of the CFIE} %

Concerning a CFIE involving the WMFIE (in the following called WCFIE), we investigate the behavior at interior resonances regarding conditioning and accuracy, with a frequency sampling step of 1\,MHz around the resonance frequencies. 
A sphere with 1\,m diameter, discretized with 999 RWG unknowns, is considered, and the reference solution is again obtained by Mie series expansion. 
The iterative solver convergence for a residual error of $10^{-4}$  is given in Fig.~\ref{fig:s2-it}. 
It is observed that the WMFIE and the WCFIE show exactly the same convergence behavior  as their classical counterparts MFIE and CFIE, where WCFIE and CFIE are implemented with a combination factor $\alpha=0.5$. Both MFIE and EFIE show worse convergence behavior at interior resonances, while the CSIE and CFIE formulations are absolutely stable. 
The advantage of the WCFIE formulation over the CSIE is clearly visible: 
Besides the reduced time per matrix vector product, the number of iterations is also reduced by a factor of two to three.
The solution accuracy is analyzed in Fig.~\ref{fig:s2-acc} by determining the maximum error in the scattered far-field. Interestingly, the WMFIE and WCFIE solutions are even more accurate than the EFIE solution for some frequencies. 

For the larger scenarios, the WCFIE is easily accelerated by the MLFMM, since only the identity operator representation is treated differently than in the standard CFIE, again with $\alpha=0.5$. 
Therefore, we examine larger scenarios in the following, starting with the 1\,m diameter sphere with $\lambda/10$ mean edge length discretization at increasing frequencies (different mesh for each simulation). 
The number of electric current unknowns ranges from $126$ to $14\,466$. 
The accuracy of the different formulations is compared in Fig.\ \ref{fig:s-acc-freq}. 
It is observed that  the large accuracy difference between the classical CFIE and the three other formulations is visible only at higher frequencies, when the shape of the sphere is discretized accurately enough.

A larger scatterer is a PEC plane Airbus A320 with $3\,025\,152$ unknowns at $300$\,MHz, with a total object length of about~$125\lambda$. 
The model is simulated for plane wave incidence with $\bm k = (0,1/\mathrm{m},-1/\mathrm{m})$ and linear $E_x$ polarization. 
For the iterative solution, an inner-outer preconditioning scheme is employed~\cite{Eibert_2007}. 
The CFIE converges with 12 iterations of the outer GMRES solver to a residual of $10^{-4}$, the WCFIE with 14 iterations, while the EFIE takes 62 iterations to converge to a residual of $3\cdot10^{-4}$ and does not converge any further within reasonable time. 
The CFIE shows a maximum relative difference of $-44.8$\,dB to the EFIE solution, while the WCFIE  is more accurate with  $-68.7$\,dB. 
Since the different far-field solutions are visually almost indistinguishable with error levels below $-40$\,dB, we also consider the CFIE with combination factor $\alpha=0.05$ instead of  $\alpha=0.5$. 
In Fig.~\ref{fig:ab}, the results for the bistatic RCS are shown for a $\varphi=90^\circ$ cut plane containing the main reflection direction at $\vartheta=135^\circ$. 
The maximum relative error of the classical CFIE-$0.05$ is at $-33.3$\,dB, the error of the new WCFIE-$0.05$ at $-56.6$\,dB, still over 20\,dB more accurate than the normal CFIE-$0.05$ and over 10\,dB more accurate than the CFIE-$0.5$ regarding maximum error. 

\begin{figure}[!t]
 \centering%
 \vspace*{-0.5cm}
 \includegraphics[]{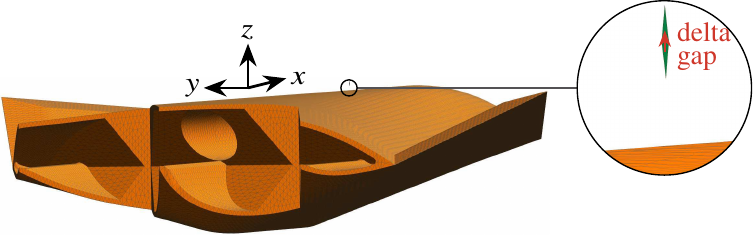}%
 \caption{Dipole voltage-gap excitation over the PEC Flamme.\label{fig:flamme}}
\end{figure}
\begin{figure}[!t]
 \centering%
 \includegraphics[]{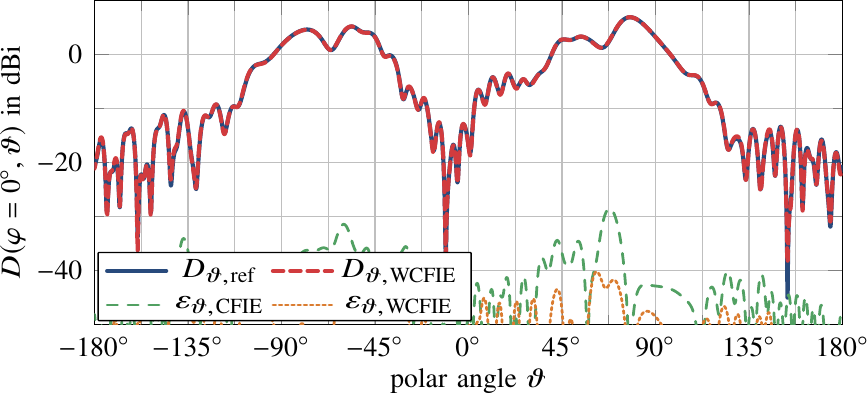}%
 \caption{Radiation characteristic of a dipole placed on top of the Flamme.\label{fig:flamme-dir}}
\end{figure}

A further large scenario is the PEC stealth object Flamme meshed with 52\,782 triangles~\cite{GuerelBagciCastelliEtAl2003,Eibert_2005}. 
The size of the structure is $16\lambda$ in $x$-direction, $6.4\lambda$ in $y$-direction and $1.6\lambda$ in $z$-direction. 
So far, only plane-wave scattering has been investigated. 
Due to the different spectral behavior, other excitations are also of interest, where in the following a dipole-like excitation on top of the structure is considered. 
It has a distance of $0.05\lambda$ in $z$-direction and it is realized by a pair of deformed triangles each with a height of $0.027\lambda$ and a width of $0.005\lambda$, see Fig.~\ref{fig:flamme}. 
The Gram matrix inversion with a diagonally preconditioned CG solver converges to a residual of $10^{-6}$ within 23 iterations on average, which takes $19$~milliseconds single-threaded on a Intel Xeon E5-1650 v4 with 3.6\,GHz. A full MLFMM matrix-vector product takes two orders of magnitude longer, about $2.2$~seconds. 
The computed directivity is given in Fig.~\ref{fig:flamme-dir}. As compared to an EFIE reference solution on a refined mesh, the maximum relative far-field errors of the normalized directivities of EFIE, CFIE, and WCFIE (both with $\alpha=0.5$) are $-46$\,dB, $-28$\,dB and $-37$\,dB, respectively.%

\section{Conclusion}

A new discretization scheme with RWG functions for the identity operator in surface integral equations was introduced based on the weak-form discretization of the combined-source condition. 
The application in the MFIE and CFIE, leading to the WMFIE and WCFIE formulations, showed a greatly enhanced accuracy in the high-frequency regime, where a further improved discretization of the $\mathcal K$ operator is not of great significance, while maintaining the good conditioning.
The main advantage of the WMFIE and the WCFIE over the RWG-discretized CSIE is the faster solution time and the utilization of Love currents, which are directly related to the tangential fields on the scatterer.

\bibliographystyle{IEEEtran}
\bibliography{IEEEabrv,ref}

\begin{IEEEbiography}
  [{\includegraphics[width=1in,height=1.25in,clip,keepaspectratio,]{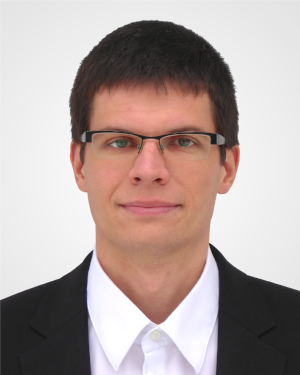}}]
  {Jonas Kornprobst} (S'17) received the B.Eng.\ degree in electrical engineering and 
  information technology from the University of Applied Sciences Rosenheim, Rosenheim, 
  Germany, in 2014 and the M.Sc.\ degree in electrical engineering and information technology from the 
  Technical University of Munich, Munich, Germany, in 2016.  Since 
  April 2016, he has been a Research Assistant with the Chair of High-Frequency 
  Engineering, Technical University of Munich.
\end{IEEEbiography}

\begin{IEEEbiography}[{\includegraphics[width=1in,height=1.25in,clip,keepaspectratio]{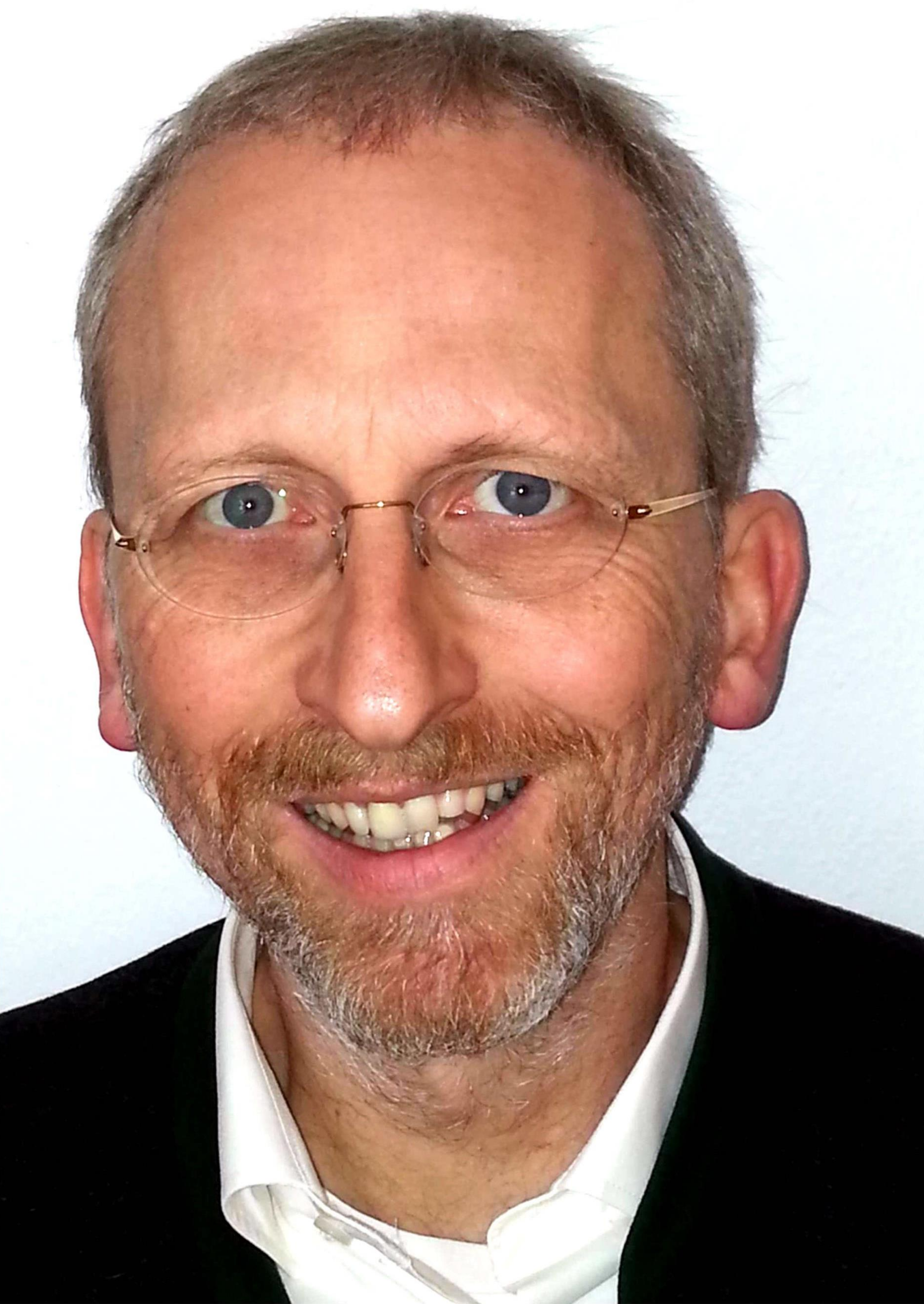}}]{Thomas F. Eibert}
             (S'93$-$M'97$-$SM'09) received the Dipl.-Ing.\,\,(FH) degree from Fachhochschule N\"urnberg, Nuremberg, Germany, the Dipl.-Ing.~degree from Ruhr-Universit\"at Bochum, Bochum, Germany, and the Dr.-Ing.~degree from Bergische Universit\"at Wuppertal, Wuppertal, Germany, in 1989, 1992, and 1997, all in electrical engineering. He is currently a Full Professor of high-frequency engineering at the Technical University of Munich, Munich, Germany.

\end{IEEEbiography}

\end{document}